\documentclass[10pt]{article}
\font \bb=msbm10 %scaled \magstep1 % Blackboard bold for real numbers
\newcommand{\R}{\hbox{\bb R}}
\newcommand{\C}{\hbox{ \bb C}}
\newcommand{\ba}{\begin{array}}
\newcommand{\ea}{\end{array}}
\newcommand{\be}{\begin{equation}}
\newcommand{\ee}{\end{equation}}
\newcommand{\bd}{\begin{displaymath}}
\newcommand{\ed}{\end{displaymath}}

\newtheorem{Theorem}{Theorem}
\newcommand{\bt}{\begin{Theorem}}
\newcommand{\et}{\end{Theorem}}
\newcommand{\br}{{\bf r}}
\begin{document}
\begin{flushleft}
\Large { \bf On Time-Dependant Symmetries
of Schr\"odinger Equation}

\end{flushleft}
\begin{flushleft}
{\it Arthur SERGHEYEV} \footnote{e-mail: arthur@apmat.freenet.kiev.ua,
arthurser@imath.kiev.ua}
\end{flushleft}
\begin{flushleft}
{\it Institute of Mathematics of the National Academy
of Sciences of Ukraine,\\
3 Tereshchenkivs'ka Street, Kyiv 4, \ Ukraine} \\

\end{flushleft}
\begin{minipage}{9cm}
\footnotesize
We show that the number of symmetry operators of order not higher that $q$
of the nonstationary $n$-dimensional ($n \leq 4$) Schr\"odinger
equation (SE) with
nonvanishing potentials is finite and does not exceed that of free
SE with zero potentials for arbitrary $q=0,1,2,\dots$. This
result is applied for the determination of the general form of time
dependance of the symmetry operators of SE with
time-independant potentials.
\end{minipage}
\section*{Introduction}
In this paper we consider (nonstationary $n$-dimensional)
Schr\"odinger equation (SE) of the form
\begin{equation} \label{se}
L \Psi \equiv (i \partial _t - \frac{1}{2} (({\bf p}-e {\bf A(\br,
t}))^{2} +  V(\br, t))) \Psi =0
\end{equation}
where $p_{a} =-i\partial /\partial x_{a} \equiv -i \partial ^{ a}$,
$a=1, \dots, n$; ${\bf p}^{2}\equiv \sum\limits_{a=1}^{n} p_{a} p_{a}$,
$\br =(x_1, \dots, x_n)$; $t, x_a \in \R$; bold letters denote
$n$-dimensional vectors.

We will study the properties of the symmetry operators of (\ref{se}),
i.e., of the linear differential operators of the form
\be \label{so}
 Q =
\sum\limits_{j=0}^{q} \lbrack \dots \lbrack F^{a_{1},\dots,a_{j}}
({\bf r},t), p_{a_1} \rbrack _{+} \dots p_{a_j} \rbrack_{+} \equiv
\sum_{j=0}^{q} b^{a_1,\dots, a_j} \partial^{a_1} \dots \partial^{a_j} ,
\ee
such that
\be \label{sym}
[L,Q]=0.
\ee
Here $[A, B]_{+} = AB+BA$, $F^{a_{1}, \dots, a_{j}}$ and $b^{a_1,
\dots, a_j}$ are symmetric with respect to the indices
$a_{1}, \dots, a_{j}$; here and below the summation over the repeated
indices of the type $a_{1}, \dots, a_{j}$ from $1$ to $n$ is understood;
the term with $j=0$ is simply $F({\bf r}, t)$.

The paper is organized as follows. In Section 1 we evaluate the
number of linearly independant time-dependant and time-independant
symmetry operators of order not higher than $q$ of SE (\ref{se}), and
in Section 2 we apply these results to establish the criterion of
existence of time-dependant symmetries for SE with time-independant
potentials. Here and below we assume $n \leq 4$.

\section{The number of symmetries}

Let us temporarily restrict ourselves to the case ${\bf A} =0$. Then
the substitution of (\ref{so}) into (\ref{sym}) yields (cf. \cite{fn})
the following equations\footnote{Dot in the above equations and below
denotes partial derivative with respect to $t$.}:
\be \label{1}
\partial^{(a_{q+1}} F^{a_{1}, \dots, a_{q})} =0,
\ee
\be \label{2}
\partial^{(a_{q}} F^{a_{1}, \dots, a_{q-1})}+2 \dot F ^{a_{1},
\dots, a_{q}} =0,
\ee
\be \label{3}
\begin{array}{lll}
\partial^{(a_{q-m}} F^{a_{1}, \dots, a_{q-m-1})}+2 \dot F ^{a_{1},
\dots, a_{q-m}}+\\
+\sum\limits_{j=0}^{\lbrack(m-1)/2 \rbrack} \frac{ 2 (-1)^{j+1}
(q-m+2j+1)!}{(q-m+1)! (2j+1)!} F^{a_{1},\dots, a_{q-m},b_{1}, \dots,
b_{2j+1}} \partial^{b_{1}} \dots \partial ^{b_{2j+1}} V=0,\\
m=2,\dots,q-1.
\end{array}
\ee
\be \label{4}
\dot F + \sum\limits_{j=0}^{\lbrack (q-1)/2 \rbrack}(-1)^{j+1}
F^{b_1,\dots, b_{2j+1}} \partial^{b_{1}} \dots \partial ^{b_{2j+1}} V
=0.
\ee
$[s]$ denotes here the integer part of the number $s$ and $(a_{1}, \dots,
a_{k})$ denotes the symmetrisation with respect to the indices
$a_{1},\dots, a_{k}$.

At first, let us show that the general
solution of the system (\ref{1})-- (\ref{4}) may not contain
arbitrary functions of $\br$.

Really, let us consider the following system\footnote{we have
introduced it by analogy with $V=0$ case \cite{fn}.}, which consists
from equation (\ref{1}) and the following differential consequences
of equations (\ref{2}), (\ref{3}):
\be \label{dc1}
 \partial^{(a_{q+1}} \partial^{a_{q}} F^{a_{1}, \dots, a_{q-1})} =0,
\ee
\be \label{dc2}
\begin{array}{lll}
\partial^{(a_{q+1}} \partial^{a_{q}} \dots \partial^{a_{j+1}} F^{a_{1},
\dots, a_{j})} =- 2 \partial^{(a_{q+1}} \partial^{a_{q}} \dots
\partial^{a_{j+2}} \dot F^{a_{1}, \dots, a_{j+1})} -\\
-\partial^{(a_{q+1}} \partial^{a_{q}} \dots \partial^{a_{j+2}}
\sum\limits_{k=0}^{\lbrack(q-j-2)/2 \rbrack} \frac{ 2 (-1)^{k+1}
(j+2k+2)!}{(j+2)! (2k+1)!} \times \\ \times F^{a_{1},\dots,
a_{j}),b_{1}, \dots, b_{2k+1}} \partial^{b_{1}} \dots \partial
^{b_{2k+1}} V,\\
j=0,\dots,q-2.
\end{array}
\ee
The general solution of the system (\ref{dc2}) (and hence, the
general solution of (\ref{1})~--~(\ref{4})) may be
represented in the form
\be \label{gs1}
F^{a_{1}, \dots, a_{j}} = F_{0}^{a_{1}, \dots, a_{j}} + G^{a_{1}, \dots,
a_{j}}, \quad j=0, \dots, q-2.
\ee
where $F_{0}^{a_{1}, \dots, a_{j}}$ is general solution of the
corresponding homogeneous equation
\be \label{hom}
\partial^{(a_{q+1}} \partial^{a_{q}} \dots \partial^{a_{j+1}}
 F^{a_{1} \dots a_{j})} =0,
\ee
i.e. the generalized Killing tensor of rank $j$ and order $q-j+1$
\cite{fn} which is (at least for $n \leq 4$) some polynomial of order
$q+1$ with respect to $\br$;
the coefficients of this polynomial are arbitrary functions of $t$
(in fact, they are not arbitrary, vide the analysis below),
and $G^{a_{1}, \dots, a_{j}}$ is some particular solution of the
corresponding inhomogeneous equation from (\ref{dc2}).
%, which vanishes if the right hand
%side of this equation equals 0 (as it is the case for $V=0$).

The representation (\ref{gs1}) is valid for $j=q-1,q$ too, but in this case
(as it follows from (\ref{1}) and (\ref{dc1})) we have
\be \label{spec}
 G^{a_{1}, \dots, a_{j}}=0,\quad j=q-1,q.
\ee
The substitution of (\ref{gs1}) (taking into account (\ref{hom})) into
(\ref{dc2}) yields the following equations:
\be \label{dc3}
\begin{array}{lll}
\partial^{(a_{q+1}} \partial^{a_{q}} \dots \partial^{a_{j+1}} G^{a_{1},
\dots, a_{j})} =- 2 \partial^{(a_{q+1}} \partial^{a_{q}} \dots
\partial^{a_{j+2}} \dot G^{a_{1}, \dots, a_{j+1})} -\\
-\partial^{(a_{q+1}} \partial^{a_{q}} \dots \partial^{a_{j+2}}
\sum\limits_{k=0}^{\lbrack(q-j-2)/2 \rbrack} \frac{ 2 (-1)^{k+1}
(j+2k+2)!}{(j+2)! (2k+1)!} \times \\ \times ( F_{0}^{a_{1},\dots,
a_{j}),b_{1}, \dots, b_{2k+1}}+ G^{a_{1},\dots, a_{j}),b_{1}, \dots,
b_{2k+1}}) \partial^{b_{1}} \dots \partial ^{b_{2k+1}} V,\\
j=0,\dots,q-2.
\end{array}
\ee
Let us require (this is obviously possible)
that the quantities $G^{a_1, \dots, a_{j}}$ must satisfy not only
(\ref{dc3}) but stronger relations
\be \label{dc3a}
\begin{array}{lll}
\partial^{(a_{j+1}} G^{a_{1},
\dots, a_{j})} =- 2 \dot G^{a_{1}, \dots, a_{j+1}} -\\
-\sum\limits_{k=0}^{\lbrack(q-j-2)/2 \rbrack} \frac{ 2 (-1)^{k+1}
(j+2k+2)!}{(j+2)! (2k+1)!} \times \\ \times ( F_{0}^{a_{1},\dots,
a_{j},b_{1}, \dots, b_{2k+1}}+ G^{a_{1},\dots, a_{j},b_{1}, \dots,
b_{2k+1}}) \partial^{b_{1}} \dots \partial ^{b_{2k+1}} V,\\
j=0,\dots,q-2,
\end{array}
\ee
and vanish if the right hand sides of (\ref{dc3a}) are equal to zero.
The comparison of (\ref{dc3a}) with (\ref{1})--(\ref{3}) yields the
following equations:
\be \label{free}
 \partial^{(a_{q-m}} F_{0}^{a_{1}, \dots, a_{q-m-1})}+2 \dot F_{0}^{a_{1},
\dots, a_{q-m}}=0, m=0, \dots, q-1
\ee
and
\be \label{extra}
\dot F_{0}= -\dot G + \sum\limits_{k=0}^{\lbrack(q-2)/2 \rbrack}
(-1)^{k} (F_{0}^{b_{1}, \dots, b_{2k+1}}+ G^{b_{1}, \dots,
b_{2k+1}}) \partial^{b_{1}} \dots \partial ^{b_{2k+1}} V.
\ee
%We see that $G^{a_{1},\dots, a_{q-2}}$ is completely determined by
%$F_{0}^{a_{1},\dots, a_{q-1},b_{1}}$ and derivatives of $V$.
Since the quantities $G^{a_1, \dots, a_{j}}$ by construction may not
contain arbitrary elements at all (except those which enter in the
quantities $F_{0}^{a_1, \dots, a_{j}}$, of course) and
the quantities $F_{0}^{a_1, \dots, a_{j}}$ may contain only arbitrary
functions of $t$, we really have
proved that the general solution of the system (\ref{1})-- (\ref{4})
may not contain arbitrary functions of $\br$.

Let us mention that the representation (\ref{gs1}) for the solution of the
system (\ref{1})~--~(\ref{4}), in which the quantities $F_{0}^{a_1,
\dots, a_{j}}$ and $G^{a_1, \dots, a_{j}}$ satisfy
(\ref{hom}), (\ref{dc3a}), (\ref{free}), (\ref{extra}), (\ref{spec}) could
be written {\it a priori} without turning to (\ref{dc2}).

The substitution of the expressions for the quantities $F_{0}^{a_1,
\dots, a_{j}}$ via generalized Killing tensors with time-dependant
coefficients into (\ref{free}) and (\ref{extra}) after equating the
coefficients at linearly independant generalized Killing tensors will
evidently yield the system of first order linear
ordinary differential equations (ODEs) with respect to $t$ for these
coefficients. Thus, these coefficients may not be arbitrary functions
of the time $t$, since the general solution of this system of ODEs may
contain only arbitrary constants at most in the same number that unknown
functions.

Thus, we have proved that the general solution of (\ref{1})--
(\ref{4}) contains neither arbitrary functions of $t$ nor arbitrary
functions of $\br$, but only arbitrary constants, and the number of
these constants does not exceed the number $\hat N_{q}^{n}$ of
coefficients of the corresponding generalized Killing tensors.
Since the system of equations (\ref{free}) and (\ref{extra}) is
linear, these arbitrary constants will enter in the quantities
$F_{0}^{a_1, \dots, a_{j}}$ (and hence in $F^{a_1, \dots, a_{j}}$) linearly.
Therefore, their number coincides with the number of linearly independant
symmetry operators of (\ref{se}).

Summing up all the above statements, we see that the number of
linearly independant symmetry operators of (\ref{se}) for arbitrary
given potential $V$ does not exceed
\be \label{nconst}
\hat N_{q}^{n} = \sum\limits_{j=0}^{q} S_{j,q}^{n}=\frac{(q+n)!
(q+n+1)!}{q! (q+1)! n! (n+1)!},
\ee
where $S_{j,q}^{n} =\frac{(j+n-1)! (q+n)!
(q-j+1)}{n! (n-1)! j! (q+1)!}$ is the number of arbirary constants in
the general time-independant solution of (\ref{hom}) \cite{fn}.

The number $\hat N_{q}^{n}$ for $n=1,2,3$ coincides with the number of
linearly independant symmetry operators of SE with $V=0, {\bf A} =0$
(and with $V= \omega \br^{2}, {\bf A} =0$ too)
of order not higher than $q$, found in \cite{fn}.
Moreover, it is straightforward to check that this coincidence takes
place for all $n=1,\dots,4$.

Now let us return to the general case of ${\bf A} \neq 0$. Reasoning
similarly to the above, we observe that the number of symmetries
of (\ref{se}) of order not higher than $q$ again does not exceed
$\hat N_q^{n}$,
since the general structure of equations for the quantities $F^{a_1,
\dots, a_{j}}$, which follow from
(\ref{sym}) is similar to (\ref{1})-- (\ref{4}): (\ref{1}) again holds true
and instead of (\ref{2})-- (\ref{4}) we have
\be \label{uzag}
\begin{array}{lll}
\partial^{(a_{q}} F^{a_{1}, \dots, a_{q-1})}+2 \dot F ^{a_{1},
\dots, a_{q}} + \dots = 0,\\
\partial^{(a_{q-m}} F^{a_{1}, \dots, a_{q-m-1})}+2 \dot F ^{a_{1},
\dots, a_{q-m}}+\dots = 0,
m=2,\dots,q-1.\\
\dot F + \dots = 0,\\
\end{array}
\ee
where dots denote some terms we need not to know explicitly; their
structure is analogous to that of (\ref{2})--(\ref{4}). In
particular, for $m$-th equation, $m=1,\dots,q$, these terms include
only $F^{a_{1}, \dots, a_{q-m+1}}, \dots F^{a_{1}, \dots, a_{q}}$ but not
$F^{a_{1}, \dots, a_{j}}$ with $j=0, \dots, q-m-1$. Therefore, we
may again represent the general solution of (\ref{uzag}) in the form
(\ref{gs1}), where now the quantities $G^{a_1, \dots, a_{j}}$ satisfy
modified version of the inhomogeneous equations (\ref{dc3a}) and
vanish if the right hand side of the corresponding equation is zero.
Repeating once more the above considerations, we obtain the following

\bt
The number of linearly independant symmetry operators of
$n$-dimensional ($n \leq 4$) SE (\ref{se}) with any {\it fixed} potentials
$V$, ${\bf A}$ of order not higher than $q$ does not exceed $\hat
N_{q}^{n}$, i.e., the number of
linearly independant symmetry operators of SE with $V=0$, ${\bf A} =0$
of order not higher than $q$.
\et

It is also interesting to evaluate the number of time-independant symmetries
of (\ref{se}) for the case when $V=V(\br)$, ${\bf A} = {\bf A}(\br)$, i.e.,
the number of time-independant operators $Q$ of the form
(\ref{so}), which commute with the Hamiltonian
\be \label{Ham}
H \equiv \frac{1}{2} (({\bf p}- e {\bf A}({\br}))^{2} + V({\bf r})).
\ee
In this case (we again temporarily set ${\bf A} =0$)
the equations (\ref{2})--(\ref{4}) read
\be \label{st}
\begin{array}{lll}
\partial^{(a_{q}} F^{a_{1}, \dots, a_{q-1})} =0,\\
\partial^{(a_{q-m}} F^{a_{1}, \dots, a_{q-m-1})}=\\
=-\sum\limits_{j=0}^{\lbrack(m-1)/2 \rbrack} \frac{2 (-1)^{j+1}
(q-m+2j+1)!}{(q-m+1)! (2j+1)!} F^{a_{1},\dots, a_{q-m},b_{1}, \dots,
b_{2j+1}} \partial^{b_{1}} \dots \partial ^{b_{2j+1}} V,\\
m=1,\dots,q-1;\\
\sum\limits_{j=0}^{\lbrack (q-1)/2 \rbrack}(-1)^{j+1} F^{a_1,\dots,
a_{2j+1}} \partial^{a_{1}} \dots \partial ^{a_{2j+1}} V =0.\\
\end{array}
\ee
Obviously, we may represent the solution of (\ref{st}) and (\ref{1}))
in the following form:
\be \label{sol}
F^{a_{1}, \dots, a_{q-l}} = \tilde F_{0}^{a_{1}, \dots, a_{q-l}} +
\tilde G^{a_{1}, \dots,
a_{q-l}}, \quad l=0, \dots, q.
\ee
where now $\tilde F_{0}^{a_{1}, \dots, a_{q-l}}$ is general solution of the
corresponding homogeneous equation
\be \label{tkk}\partial ^{(a_{q-l+1}} F^{a_{1} \dots a_{q-l})} =0,
\ee
i.e. the generalized Killing tensor of rank $q-l$ and order 1 \cite{fn}
which is (for $n \leq 4$) some polynomial of order $q-l$ with
respect to $\br$, containing $K_{j}^{n} = \frac{(j+n-1)! (j+n)!}{j! (j+1)!
(n-1)! n!}$ arbitrary constants \cite{fn},
and $\tilde G^{a_{1}, \dots, a_{q-l}}$ is particular solution of the
corresponding inhomogeneous equation from (\ref{st}) ($\tilde G^{a_{1},
\dots, a_{j}}=0$, $j=q-1,q$ in virtue of (\ref{1}) and of the first
line of (\ref{st})).

In complete analogy with the above it is clear that the general
solution of the system (\ref{st}) and (\ref{1}) may depend at most from
\be \label{nst}
\tilde N_{q}^{n} = \sum\limits_{j=0}^{q} K_{j}^{n},
\ee
arbitrary constants, which again enter in it linearly; it really
depends on $\tilde N_{q}^{n}$ constants for, e.g., $V=0, {\bf A}=0$.

The explicit calculations show that
\begin{equation}
\begin{array}{lll}
\tilde N_{q}^{1} = q + 1,
\tilde N_{q}^{n} = \frac{(q+n+1)!}{q! (2n-1)!} P_{n} (q), n =2, 3, 4,\\
\mbox{where} \quad P_{2}(q)=1, P_{3}(q) = 2 q + 5, P_{4}(q) =5 q^{2}
+ 30 q +42.\\
\end{array}
\end{equation}

In complete analogy with above considerations, we may return to the
general case, when ${\bf A} \neq 0$. Thus, the following statement
holds true:
\bt
For any {\it fixed} $V(\br)$ and ${\bf A}(\br)$ the number of
time-independant operators of the form (\ref{so}) of order not higher
than $q$, which commute with $H$ (\ref{Ham}), does not exceed $\tilde
N_{q}^{n}$, i.e., the same number for $V=0, {\bf A} =0$, for $n=1,
\dots, 4$.
\et

If we suppose that the quantities $F^{a_1, \dots, a_{j}}$, ${\bf A},
V$ are "generalized" distributions\footnote{the products and
derivatives of these "generalized" distributions again are some
"generalized" distributions.}
(e.g., in Co\-lom\-beau's sense) \cite{eg}, both our theorems hold
true, since all the necessary derivatives exist (of course, as
"generalized" distributions) and hence we may apply our
results to even such singular potentials as Dirac delta-functions and
other distributions.

 Let us also mention that our results exploit {\it ad hoc}
adjusted ideas of Cartan's theory of compatibility of overdetermined
systems of partial differential equations \cite{cart} for the system
(\ref{1}) and (\ref{uzag}).

\section{Time-dependant Symmetries}
In this section we consider the case when both $V$ and $\bf A$ are
time-independant.

The symmetry operators of the form (\ref{so}) of (\ref{se}) may be
considered as the elements of Lie algebra $Sym$ of all the generalized
symmetries of SE (cf.\cite{o}), since the Lie algebra (with
respect to the usual commutator $[A, B] = A B - B A$) $Lin$ of all the
operators of the form (\ref{so}) may be homomorphically mapped (with
zero kernel) to the Lie algebra $S$ (with respect to Lie bracket) of
all the differential operators of the form
\be
{\rm K} = \eta (\br,t,\Psi,\dots, \Psi^{(q)}) \partial /
\partial \Psi
\ee
by mapping $Q$ (\ref{so}) into
\begin{equation} \label{homo}
Q \rightarrow {\rm Q}'= (\sum_{j=0}^{q}
b^{a_1,\dots, a_j} \partial^{a_1} \dots \partial^{a_j} \Psi)
\partial/\partial \Psi
\end{equation}
Here $\Psi ^{(q)}$ denotes the set of derivatives of $\Psi$ with
respect to $\br$ of order $q$.

Let $W$ be the image of $Lin$ under (\ref{homo}).
% Moreover, in what follows
%we will simply identify $W$ with $Lin$ for the sake of simplicity.
 For such a $W$ the dimensions $v^{(q)}$
of the spaces of symmetry operators of order not higher than $q$
$V^{(q)} = W \bigcap Sym^{(q)}$ (as in \cite{as}, we set
$V=W \bigcap Sym$) are finite and do not exceed $\hat N_{q}^{n}$, as
we have established in the Theorem 1 above. Since in the case
considered (\ref{se}) admits Lie symmetry $\partial/\partial t$ and
if $R$ is symmetry operator of (\ref{se}), then so does $\partial
R/\partial t$, all the conditions of Theorem 1 from \cite{as} are
fulfilled.
Therefore, we may seek for all the linearly independant symmetry
operators of (\ref{se}) of order $q$ in the form:
\begin{equation} \label{ansatz}
R =\exp (\lambda t) \sum_{k=0}^{m} C_{k} \frac{t^{k}}{k!}
\end{equation}
where $\lambda \in \C$ and $C_k$ are some differential operators from $Lin$
with time-in\-de\-pendant coefficients, $m \leq v^{(q)}-1$.

Moreover, applying Theorem 2 from \cite{as2}, we immediately obtain
the following result, which holds true for $n \leq 4$:
\bt
SE (\ref{se}) with time-independant potentials $V$ and $\bf A$ possesses
time-dependant symmetry operators if and only if it possesses (at least one)
symmetry operator of the form
\be \label{case1}
Q= \exp(\lambda t) K_{0}, \lambda \in \C, \lambda \neq 0
\ee
or of the form
\be \label{case2}
 Q=  K_{0} + t K_{1},\:  K_1 \neq 0,
\ee
where $K_0, K_1$ are differential operators from $Lin$ with time-independant
coefficients.
\et

Now let us write down the commutation relations for $C_k$,
obtained as a result of substitution of (\ref{ansatz})
into (\ref{sym}):
\be \label{cr1}
[H, C_m]= i\lambda C_m; \quad [H, C_l] = i \lambda C_l +i C_{l+1},
l=0, \dots, m-1.
\ee
Having exluded $C_l$, $l=1, \dots, m$, we may rewrite (\ref{cr1})
simply as
\be \label{cr2}
 (-i {\rm ad}_H - \lambda)^{m} C_0 =0,
\ee
where ${\rm ad}_H R \equiv [H,R]$.
Similarly, for $K_0$ and $K_1$ we have
\begin{eqnarray} \label{cr3}
[ H, K_{0}] = i\lambda K_{0} \; \mbox{for the case (\ref{case1})}
\\[\medskipamount]
\label{cr5}
 [ H, K_0 ] =i K_1, \: [ H, K_1] =0 \; \mbox{for the case (\ref{case2})}.
\end{eqnarray}
For the case (\ref{case2}) from (\ref{cr5}) it follows that
\be \label{cr4}
[H, [H, K_{0}]]=0,
\ee
and hence by analogy with the theory of mastersymmetries of
integrable equations \cite{fu} it is natural to call $K_0$ the
\textit{mastersymmetry}
of Schr\"odinger equation~(\ref{se}).

\end{document}